\documentclass[12pt,leqno]{article}
\usepackage{amsfonts}
\usepackage{amsmath, amsthm, amsfonts, amssymb, color}
\usepackage{mathrsfs}
\usepackage{color}
\usepackage[notcite,notref]{showkeys}

\theoremstyle{definition}

\newcommand{\scr}[1]{\mathscr #1}
\definecolor{wco}{rgb}{0.5,0.2,0.3}

\numberwithin{equation}{section} \theoremstyle{remark}

\newcommand{\ua}{\uparrow}

\title{{\bf  Functional SPDE    with Multiplicative Noise and Dini  Drift }\footnote{Supported in
 part by  NNSFC(11131003, 11431014), the 985 project and the Laboratory of Mathematical and  Complex Systems.} }
\author{
{\bf     Feng-Yu Wang $^{a), b)}$, Xing Huang $^{a)}$  }\\
\footnotesize{ a)School of Mathematical Sciences,
Beijing Normal
University, Beijing 100875, China}\\
 \footnotesize{ b)Department of Mathematics,
Swansea University, Singleton Park, SA2 8PP, United Kingdom}\\
\footnotesize{  wangfy@bnu.edu.cn, F.-Y.Wang@swansea.ac.uk}}
\begin{document}
\allowdisplaybreaks
\def\R{\mathbb R}  \def\ff{\frac} \def\ss{\sqrt} \def\B{\mathbf
B}
\def\N{\mathbb N} \def\kk{\kappa} \def\m{{\bf m}}
\def\ee{\varepsilon}\def\ddd{D^*}
\def\dd{\delta} \def\DD{\Delta} \def\vv{\varepsilon} \def\rr{\rho}
\def\<{\langle} \def\>{\rangle} \def\GG{\Gamma} \def\gg{\gamma}
  \def\nn{\nabla} \def\pp{\partial} \def\E{\mathbb E}
\def\d{\text{\rm{d}}} \def\bb{\beta} \def\aa{\alpha} \def\D{\scr D}
  \def\si{\sigma} \def\ess{\text{\rm{ess}}}
\def\beg{\begin} \def\beq{\begin{equation}}  \def\F{\scr F}
\def\Ric{\text{\rm{Ric}}} \def\Hess{\text{\rm{Hess}}}
\def\e{\text{\rm{e}}} \def\ua{\underline a} \def\OO{\Omega}  \def\oo{\omega}
 \def\tt{\tilde} \def\Ric{\text{\rm{Ric}}}
\def\cut{\text{\rm{cut}}} \def\P{\mathbb P} \def\ifn{I_n(f^{\bigotimes n})}
\def\C{\scr C}      \def\aaa{\mathbf{r}}     \def\r{r}
\def\gap{\text{\rm{gap}}} \def\prr{\pi_{{\bf m},\varrho}}  \def\r{\mathbf r}
\def\Z{\mathbb Z} \def\vrr{\varrho} \def\ll{\lambda}
\def\L{\scr L}\def\Tt{\tt} \def\TT{\tt}\def\II{\mathbb I}
\def\i{{\rm in}}\def\Sect{{\rm Sect}}  \def\H{\mathbb H}
\def\M{\scr M}\def\Q{\mathbb Q} \def\texto{\text{o}} \def\LL{\Lambda}
\def\Rank{{\rm Rank}} \def\B{\scr B} \def\i{{\rm i}} \def\HR{\hat{\R}^d}
\def\to{\rightarrow}\def\l{\ell}\def\iint{\int}
\def\EE{\scr E}\def\no{\nonumber}
\def\A{\scr A}\def\V{\mathbb V}\def\osc{{\rm osc}}
\def\BB{\scr B}\def\Ent{{\rm Ent}}
\def\U{\scr U}

\maketitle

\begin{abstract} Existence, uniqueness and non-explosion of the mild solution are proved for  a class of   semi-linear functional SPDEs with  multiplicative noise and   Dini continuous drifts. In the finite-dimensional and  bounded time delay setting,  the log-Harnack inequality and $L^2$-gradient estimate are derived.   As the Markov semigroup is associated to the functional (segment) solution of the equation,    one needs to make analysis  on the  path space  of the solution in the time interval of delay.  
\end{abstract} \noindent
 \vskip 2cm

\section{Introduction}

It is well known by Dominique Bakry and his collaborators that the curvature lower bound condition of a diffusion process is equivalent to a number of  gradient inequalities for the associated Markov semigroup, see e.g. the recent monographs  \cite{B, Wbook}. Among many other equivalent inequalities, the $L^2$-gradient estimate of type
$$|\nn P_t f|^2 \le C(t) P_t f^2,\ \ t>0$$ has been extended to more general situations without curvature conditions, see e.g. \cite{PW, RW, SWY, W10, WZ} and references within for the study of   SDEs/SPDEs with non-Lipschitz coefficients. The $L^2$-gradient estimate links to    the log-Harnack inequality which has further applications in analysis of Markov operators, see \cite{WB} and references within. Recently, by constructing  a Zvonkin  type transformation in Hilbert spaces, the $L^2$-gradient estimate and log-Harnack inequality have been derived in \cite{W} for semi-linear SPDEs with Dini drifts. In the present paper we aim to extend these results to SPDEs with delay.

We will consider semi-linear SPDEs with delay in a separable Hilbert space $(\mathbb{H}, \<\cdot,\cdot\>, |\cdot|)$.  To describe the time delay, let $\nu$ be a non-trivial  measure on $(-\infty,0)$ such that
\beq\label{nu} \nu \ \text{is\ locally\ finite\ and}\   \nu(\cdot-t)\le \kk(t)\nu(\cdot), \ \  t> 0\end{equation}   for some  increasing  function $\kk: (0,\infty)\to (0,\infty)$.
This condition is crucial to prove the pathwise (see the proof of Proposition  \ref{P3.3} below), and to determine  the state space of the segment solutions (see Remark 1.1  below).  Obviously, \eqref{nu} holds for $\nu(\d\theta):=1_{(-\infty,0)}(\theta)\rr(\theta)\d\theta$ with density $\rr\ge 0$ satisfying
$\rr(\theta-t)\le\kk(t)\rr(\theta)$ for $\theta<0,$ which is the case if, for instance, $\rr(\theta)= \e^{\ll\theta} 1_{[-r_0,0)}(\theta)$ for some constants
$\ll\in\R$ and $r_0\in (0,\infty].$ 
Then the  state space of the segment process under study  is given by
$$\C_\nu:= \bigg\{\xi: (-\infty,0]\to \H \ \text{is\ measurable\ with}\ \nu(|\xi|^2)<\infty\bigg\},$$
 where $\nu(f):=\int_{-\infty}^0 f(\theta)\nu(\d\theta)$ for $f\in L^1(\nu).$
 Let
$$\|\xi\|_{\C_\nu}=\ss{ \nu (|\xi|^2) + |\xi(0)|^2},\ \ \xi\in \C_\nu.$$ Throughout the paper, we   identify  $\xi$ and $\eta$ in $\C_\nu$ if   $\xi=\eta\ \nu$-a.e. and $\xi(0)=\eta(0)$, so that
$\C_\nu$ is a separable Hilbert space with inner produce
$$\<\xi,\eta\>_{\C_\nu}:= \nu(\<\xi,\eta\>) +\xi(0)\eta(0),\ \ \xi, \eta\in \C_\nu.$$

For a map $X: \R\to \H$ and $t\ge 0$, let $X_t: (-\infty, 0]\to \H$ be defined by
$$X_t(\theta)= X(t+\theta),\ \ \theta\in (-\infty,0],$$ which describes the path of $X$ from $-\infty$ to time $t$.  We call $X_t$ the segment of $X$ at time $t$.

Consider the following semi-linear SPDE on $\H$:
\beq\label{E1} \d X(t)= \big\{A X(t)+b(t,X(t))+B(t,X_{t})\big\}\d t+Q(t,X(t))\d W(t),
\end{equation}
where \beg{enumerate} \item[$\odot$] $(A,\D(A))$ is a negative definite self-adjoint operator on $\mathbb{H}$;
\item[$\odot$]  $B: [0,\infty)\times \C_\nu\to \mathbb{H}$ and  $b: [0,\infty)\times \mathbb{H}\to \mathbb{H}$ are measurable and locally bounded (i.e. bounded on bounded sets);
\item[$\odot$] $W=(W(t))_{t\geq 0}$ is a cylindrical Brownian motion on a separable Hilbert space $\mathbb{\bar{H}}$, with respect to a complete filtration  probability space $(\OO, \F, \{\F_{t}\}_{t\ge 0}, \P)$. More precisely, $W(t):=\sum_{n=1}^{\infty}{B^{n}(t)\bar{e}_{n}}$ for a sequence of independent one-dimensional Brownian motions $\{B^{n}(t)\}_{n\geq 1}$ with respect to $(\OO, \F,
\{\F_{t}\}_{t\ge 0}, \P)$, and   an orthonormal basis $\{\bar{e}_{n}\}_{n\geq 1}$ on $\mathbb{\bar{H}}$;
\item[$\odot$]  $Q: [0,\infty)\times \mathbb{H}\to \L(\mathbb{\bar{H}}; \mathbb{H})$ is measurable, where $\L(\mathbb{\bar{H}};\mathbb{H})$ is the space of bounded linear operators from $\mathbb{\bar{H}}$ to $\mathbb{H}$. \end{enumerate}

\beg{defn} For any $\xi\in \C_\nu$,  an adapted continuous process $(X(t))_{t\in [0,\zeta)}$ on $\H$ is called a mild solution to \eqref{E1} with initial value $X_0=\xi$ and life time $\zeta$,
if $\zeta$ is a stopping time such that
$\limsup_{t\uparrow \zeta} |X(t)|=\infty$ holds on $\{\zeta<\infty\}$, the Lebesgue integral $\int_0^t \e^{(t-s)A} \big\{b(s, X(s))+ B(s, X_s)\}\d s$ and  the It\^o integral $ \int_0^t \e^{(t-s)A} Q(s, X(s))\d W(s)$ are well defined on $\H$ for $t\in [0, \zeta)$, and
\beq\label{E2} \beg{split} X(t)=\ &\e^{At} \xi(0) +\int_0^t \e^{(t-s)A} \big\{b(s, X(s))
 + B(s, X_s)\}\d s \\ &\qquad + \int_0^t \e^{(t-s)A} Q(s, X(s))\d W(s),\ \ t\in [0, \zeta)\end{split}\end{equation} holds. Here, due to $X_0=\xi$, 
  $X$ is extended to $(-\infty,0)$  with $X(\theta)= \xi(\theta)$ for $\theta\le 0$.   If the solution exists uniquely, we denote it by
$(X^\xi(t))_{t\in [0,\zeta)}.$ The solution is called non-explosive if $\zeta=\infty$ a.s.   \end{defn}

\paragraph{Remark 1.1} We note that condition \eqref{nu} ensures that
    the segment solution $(X_t^\xi)_{t\in [0,\zeta)}$ is an adapted   process on $\C_\nu$.  If, for any initial value  the equation \eqref{E1} has a unique non-explosive mild solution, then let
$$P_t f(\xi)= \E f(X_t^\xi),\ \ t\ge 0, f\in \B_b(\C_\nu),$$ where $\B_b(\C_\nu)$ is the set of all bounded measurable functions on $\C_\nu$.
 In the time-homogenous case (i.e. $b(s,\cdot), B(s,\cdot)$ and $Q(s,\cdot)$ do not depend on $s$),  $P_t$ is a Markov semigroup on $\B_b(\C_\nu)$. In general, for any $s\ge 0$,
 let $X^\xi_s(t)$ be the mild solution of the equation \eqref{E1}   for $t\ge s$ with $X_s=\xi$,  then the associated     Markov semigroup  $\{P_{s,t}\}_{t\ge s\ge 0}$  is defined by
 $$P_{s,t}f(\xi) = \E f(X_{s,t}^\xi),\ \ t\ge s, f\in \B_b(\C_\nu),$$ where $X^\xi_{s,t}$ is the segment of $X_s^\xi(\cdot)$ at time $t$.

\

The remainder of the paper is organized as follows. In Section 2 we study the existence, uniqueness and non-explosion of the mild solution.   In Section 3 we investigate  the
log-Harnack inequality and $L^2$-gradient estimate for $P_{s,t}$ when $\H$ is finite-dimensional and supp$\,\nu\subset [-r_0,0]$ for some constant $r_0\in (0,\infty).$ Explanations on making  these restrictions are given in the beginning of Section 3. 

  \section{Existence, uniqueness and  non-explosion}

Let $\|\cdot\|$ and $\|\cdot\|_{HS}$ denote, respectively,  the operator norm and the Hilbert-Schmidt norm for linear operators, and let $\L_{HS}(\mathbb{\bar{H}}; \mathbb{H})$ be the space of   Hilbert-Schmidt linear operators from $\mathbb{\bar{H}}$ to $\mathbb{H}$. Moreover, let
\beg{equation*}\beg{split}
\D= \Big\{\phi: [0,\infty)\to [0,\infty) \text{ is increasing}, \phi^{2} \text{ is concave}, \int_0^1{\frac{\phi(s)}{s}\d s}<\infty\Big\}.
\end{split}\end{equation*}
As in \cite{W}, we will use this class of functions  to characterize the Dini continuity of the drift $b$. Note that  the condition $\int_0^1{\frac{\phi(s)}{s}\d s}<\infty$ is   known as Dini condition, due to the notion of Dini continuity.

\

 To ensure the existence and uniqueness of solutions, we make the following assumptions:
\beg{enumerate}
\item[$(A1)$] There exists $\vv\in (0,1)$ such that   $(-A)^{\varepsilon-1}$ is of trace class for some $\varepsilon \in(0,1)$; i.e. $\sum_{n=1}^{\infty}{\lambda_{n}^{\varepsilon-1}}<\infty$ for $0< \lambda_{1}\leq \lambda_{2}\leq\cdots\cdots$ being all eigenvalues of $-A$ counting multiplicities. Let $\{e_n\}_{n\ge 1}$ be the eigenbasis of $A$ associated with   $\{\ll_n\}_{n\ge 1}. $
\item[$(A2)$]  $Q(t,\cdot)\in C^{2}(\mathbb{H};\L(\mathbb{\bar{H}};\mathbb{H}))$ for  $t\in[0,\infty),  \  Q(t,x) Q(t,x)^{\ast}$ is invertible for   $(t,x)\in[0,\infty)\times \mathbb{H}$, and
\beg{equation*}\beg{split}
\|\nabla Q(t,x)\|+\|\nabla^{2} Q(t,x)\|+\| Q(t,x)\|+\|\{Q(t,x) Q(t,x)^{\ast}\}^{-1}\|
\end{split}\end{equation*}is locally bounded in $(t,x)\in[0,\infty)\times \mathbb{H}$, where $\nn$ is the gradient operator on $\H$.
Moreover,  for any $x\in\mathbb{H}$ and $t\geq 0$,
\beq\label{1.3}
\lim_{n\to \infty}\|Q(t,x)-Q(t,\pi_{n}x)\|_{HS}^{2}:=\lim_{n\to \infty}\sum_{k\geq1}|(Q(t,x)-Q(t,\pi_{n}x))\bar{e}_{k}|^{2}=0,
\end{equation} where $\pi_n: \H\to \H_n:=\text{span}\{e_i: 1\le i\le n\}$ are orthogonal projections for $n\ge 1.$
\item[$(A3)$] For any $n\geq1$, there exits $\phi_{n}\in\D$ such that
\beq\label{1.2}
|b(t,x)-b(t,y)|\leq \phi_{n}(|x-y|),\quad t\in[0,n], x,y\in \mathbb{H},\ \text{with}\ |x|\vee |y|\leq n.
\end{equation}
\item[$(A4)$] For any $n\ge 1$ there exists a constant $C_n\in (0,\infty)$ such that
$$|B(t,\xi)-B(t,\eta)|^2\le C_n\|\xi-\eta\|_{\C_\nu}^2,\ \ t\in [0,n], \xi,\eta\in \C_\nu\ \text{with}\   \|\xi\|_{\C_\nu}\lor \|\eta\|_{\C_\nu}\le n.$$
\end{enumerate}

When the delay term $B$ vanishes,   the existence and uniqueness of mild solutions  have  been proved in \cite{W} under assumptions $(A1)$-$(A3)$. The
additional assumption $(A4)$ means that the delay term is  locally Lipschitzian in $\C_\nu$.      Note that this condition allows unbounded   time delay, i.e. supp$\,\nu$ might be unbounded.

Fo $T>0$, let $\|\cdot\|_{T,\infty}$ denote the uniform norm on $[0,T]\times\H.$   The main result of this section is the following.

\beg{thm}\label{T2.1} Assume that $\eqref{nu}$ and $(A1)${\rm -}$(A4)$ hold. Then:
\beg{enumerate}
\item[$(1)$]   For any initial value $X_{0}\in  \C_{\nu}$, the equation $\eqref{E1}$ has unique mild solution $(X(t))_{t\in[0,\zeta)}$ with life time $\zeta$.
\item[$(2)$]  Let $\|Q\|_{T,\infty}<\infty$ for $T\in (0,\infty)$. If  there exist two positive increasing functions $\Phi, h:[0,\infty)\times[0,\infty)\to(0,\infty)$ such that
 $\int_{1}^{\infty}\frac{\d s}{\Phi_t(s)}=\infty$ for $t>0$ and
\beq\label{2.1}
\langle B(t,\xi+\eta)+b(t,(\xi+\eta)(0)),\xi(0)\rangle \leq\Phi_{t}(\|\xi\|_{\C_\nu}^{2})+h_{t}(\|\eta\|_{\C_\nu}), \xi, \eta\in\C_\nu, t\geq 0,
\end{equation}
then the mild solution is non-explosive.
\end{enumerate}
\end{thm}

In Subsection 2.1 we investigate the pathwsie uniqueness, which, together with the weak existence, implies the existence and uniqueness of mild solutions according to  the Yamada-Watanabe principle.
Complete proof of Theorem \ref{T2.1} is addressed in Subsection 2.2.

\subsection{Pathwise uniqueness}

In this section, we prove the pathwise uniqueness of the mild solution to \eqref{E1} under   $(A1)$ and the following stronger versions of $(A2)$-$(A4)$:
\beg{enumerate}
\item[$(A2')$] In addition to $(A2)$ there holds
$$\|\nabla Q \|_{T,\infty}+\|\nabla^{2} Q\|_{T,\infty}+\|Q\|_{T,\infty}+\|(Q Q^{\ast})^{-1}\|_{T,\infty}<\infty,\ \ T>0.$$
\item[$(A3')$] For any $T>0$,  $\|b\|_{T,\infty}<\infty$   and   there exits $\phi\in\D$ such that
$$
|b(t,x)-b(t,y)|\leq \phi(|x-y|),\ \ t\in[0,T], x,y\in \mathbb{H}.
$$
\item[$(A4')$] For any $T>0$ there exists a constant $C \in (0,\infty)$ such that
$$|B(t,\xi)-B(t,\eta)|^2\le C \|\xi-\eta\|_{\C_\nu}^2,\ \ t\in [0,T], \xi,\eta\in \C_\nu.$$
 \end{enumerate}

\beg{prp}\label{P3.3} Assume $\eqref{nu}$, $(A1)$ and $(A2')$-$(A4')$.   Let $X(t)_{t\geq 0}$ and $Y(t)_{t\geq 0}$ be two adapted continuous processes on $\mathbb{H}$ with $X_{0}=Y_{0}=\xi\in\C_\nu$. For any $n\geq 1$, let
\beg{equation*}
\tau_{n}^{X}=n\wedge \inf\{t\geq 0:|X(t)|\geq n\},\ \ \tau_{n}^{Y}=n\wedge \inf\{t\geq 0:|Y(t)|\geq n\}.
\end{equation*}
If  for all $t\in[0,\tau_{n}^{X}\wedge\tau_{n}^{Y}]$
\beg{equation}\label{*1}\beg{split}
&X(t)= \e^{At} \xi(0)+\int_{0}^{t} \e^{A(t-s)}(b(s,X(s))+B(s,X_{s}))\d s+\int_{0}^{t} \e^{A(t-s)}Q(s,X(s))\d W(s),\\
&Y(t)= \e^{At} \xi(0)+\int_{0}^{t} \e^{A(t-s)}(b(s,Y(s))+B(s,Y_{s}))\d s+\int_{0}^{t} \e^{A(t-s)}Q(s,Y(s))\d W(s),
\end{split}\end{equation}
then   $X(t)=Y(t)$  for all $t\in[0,\tau_{n}^{X}\wedge\tau_{n}^{Y}]$. In particular,   $\tau_{n}^{X}=\tau_{n}^{Y}$.
\end{prp}

We will prove this result by using the Zvonkin type  transform constructed in \cite{W}. Let $\{P_{s,t}^0\}_{t\ge s\ge 0}$ be the Markov semigroup associated to the O-U type SDE on $\H$:
$$\d Z_{s}(t)= AZ_{s}(t)\d t + Q(t, Z_{s}(t))\d W(t),\ \ t\ge s.$$ Given $\ll,T>0$, consider the following equation on $\H$:
\beq\label{E2} u(s,\cdot)= \int_s^T\e^{-\ll(t-s)}P_{s,t}^0 \big\{\nn_{b(t,\cdot)}u(t,\cdot) +b(t,\cdot)\big\}\d t,\ \ s\in [0,T].\end{equation}
The next   result is essentially due to \cite{W}, where the first assertion follows from \cite[Lemma 2.3]{W} and the second can be proved as in the proof of \cite[Proposition 2.5]{W} by taking into account the additional drift $B$. So, to save space we skip the proof.

\beg{lem}[\cite{W}]\label{L3.1}  Assume $(A1)$ and $(A2')$-$(A4')$.  For any $T>0$,  there exists a constant $\lambda(T)>0$ such that 
the following assertions hold for  $\ll\ge \ll(T)$: 
\beg{enumerate} \item[$(1)$]  $\eqref{E2}$ has a unique solution $u\in C([0,T]; C_{b}^{1}(\mathbb{H}; \mathbb{H}))$ and
$$
\lim_{\lambda\to \infty}\big\{\|u\|_{T,\infty}+\|\nabla u\|_{T,\infty}+\|\nabla^{2} u\|_{T,\infty}\big\}=0.
$$
\item[$(2)$]  Let $\tau$ be a stopping time. If an adapted continuous process $(X(t))_{t\in[0,T\wedge\tau]}$ on $\mathbb{H}$ satisfies
$$
X(t)= \e^{At} X(0)+\int_{0}^{t} \e^{A(t-s)}(b(s,X(s))+B(s,X_{s}))\d s
+\int_{0}^{t} \e^{A(t-s)}Q(s,X(s))\d W(s)$$  for all $t\in[0,\tau\wedge T],$
then  
\beg{equation*} \beg{split}
X(t)=\ &\e^{At}\big\{X(0)+u(0,X(0)))-u(t,X(t)\big\}\\
&+\int_{0}^{t} \e^{A(t-s)}\big\{Q(s,\cdot)
 +(\nabla u(s,\cdot))Q(s,\cdot)\big\}(X(s))\d W(s)\\
 &+\int_{0}^{t}\Big\{(\lambda-A)\e^{A(t-s)}u(s,X(s))
 +\e^{A(t-s)}\big[B(s,X_{s})+\nabla_{B(s,X_{s})}u(s,X(s))\big]\Big\}\d s
\end{split}\end{equation*} holds for all $t\in [0,\tau\land T]. $\end{enumerate}
\end{lem}

\beg{proof}[Proof of Proposition $\ref{P3.3}$] For any $m\geq 1$, let $\tau_{m}=\tau_{n}^{X}\wedge\tau_{n}^{Y}\wedge \inf\{t\geq0: |X(t)-Y(t)|\geq m\}$. It suffices to prove that for any $T>0$ and $m\geq1$,
\beq\label{3.5}
\int_{0}^{T} \mathbb{E}\big\{1_{\{s<\tau_{m}\}}|X(s)-Y(s)|^{2}\big\}\d s=0.
\end{equation}
Let $\lambda>0$ be such that assertions in Lemma \ref{L3.1}  hold. Let $I$ be the identity operator. Due to \eqref{*1}, Lemma \ref{L3.1}(2)  with $\tau=\tau_{m}$ implies
\beq\label{X1}X(t)-Y(t)= \LL(t)+\Xi(t),\ \ t\in [0,\tau_m\land T],\end{equation}
where
\beg{equation*} \beg{split}\LL(t)&:= \int_0^t \e^{(t-s)A} \Big\{\big(I+\nn u(s, X(s))\big)\big(B(s, X_s)- B(s,Y_s)\big) \\
&\qquad\qquad\qquad\qquad +\big(\nn u(s, X(s))-\nn u(s, Y(s))\big)B(s,Y_s)\Big\}\d s,\\
\Xi(t)&:= u(t,Y(t))-u(t,X(t))+\int_{0}^{t}(\lambda-A)\e^{A(t-s)}\big\{u(s,X(s))-u(s,Y(s))\big\}\d s\\
&\quad +\int_{0}^{t} \e^{A(t-s)}\big\{\nabla u(s,X(s))-\nabla u(s,Y(s))\big\}Q(s,X(s))\d W(s)\\
&\quad +\int_{0}^{t} \e^{A(t-s)}\big(\nabla u(s,Y(s))+I\big)\big(Q(s,X(s))-Q(s,Y(s))\big)\d W(s),t\in[0,\tau_{m}\wedge T].\end{split}\end{equation*}
According to the proof of \cite[Proposition 3.1]{W} (see the inequality before (3.8) therein), when $\ll\ge\ll(T)$ is large enough there exists a constant $C_0\in (0,\infty)$ such that
\beq\label{X2}  \int_0^r \e^{-2\ll t} \E\big[1_{\{t<\tau_m\}} |\Xi(t)|^2\big]\,\d t\le \ff 3 4 \GG(r) + C_0 \int_0^r \GG(t)\d t,\ \ r\in [0,T]\end{equation} holds for
\beq\label{X3} \GG(t):= \int_0^t \e^{-2\ll s}  \E\big[1_{\{s<\tau_m\}} |X(s)-Y(s)|^2\big]\,\d s,\ \ t\in [0,T],\end{equation} which is denoted by $\eta_t$ in \cite{W}.
So, to prove \eqref{3.5}, it remains to estimate the corresponding term for $\LL(t)$ in place of $\Xi(t).$
Noting that $X_0=Y_0$ in $\C_\nu$ implies $X=Y\ \nu$-a.e. on $(-\infty,0)$, by \eqref{nu} we have
$$\int_{-\infty}^{-s} |X(s+q)-Y(s+q)|^2 \nu(\d q) =\int_{-\infty}^0 |X(\theta)-Y(\theta)|^2 \nu(\d\theta-s)=0,\ \ s\ge 0.$$ 
So, by $\|\e^{A(t-s)}\|\le 1$ for $t\ge s$, Lemma \ref{L3.1}(1)  and  $(A4')$, we may find   constants $C_1, C_2\in (0,\infty)$ such that
  \beg{equation*}\beg{split} |\LL(t)|^2 &
 \leq C_1\int_{0}^{t}\big\{|B(s,X_{s})-B(s,Y_{s})|^{2}+|X(s)-Y(s)|^2\big\} \d s\\
&\leq C_2 \int_{0}^{t}|X(s)-Y(s)|^{2}\d s+C_2\int_{0}^{t}\d s\int_{-\infty}^{0}|X(s+q)-Y(s+q)|^{2}\nu(\d q) \\
&=  C_2 \int_{0}^{t}|X(s)-Y(s)|^{2}\d s+C_2\int_{0}^{t}\d s\int_{-s}^{0}|X(s+q)-Y(s+q)|^{2}\nu(\d q) \\
&=  C_2 \int_{0}^{t}|X(s)-Y(s)|^{2}\d s+C_2\int_{-t}^{0}\nu(\d q)\int_{-q}^t|X(s+q)-Y(s+q)|^{2}\d s \\
&\le K(T) \int_0^t |X(s)-Y(s)|^2\d s,\ \ t\in [0,T],
\end{split}\end{equation*} where $K(T):= C_2+C_2 \nu([-T,0))<\infty$ since $\nu$ is locally finite by \eqref{nu}.  Thus,
\beg{equation*}\beg{split} \int_0^r \e^{-2\ll t} \E\big[1_{\{t<\tau_m\}} |\LL(t)|^2\big]\,\d t&\le K(T) \E \int_0^r \e^{-2\ll t}1_{\{t<\tau_m\}}\d t \int_0^t |X(s)-Y(s)|^2\d s\\
&\le K(T)\int_0^r \GG(t)\d t,\ \ r\in [0,T].\end{split}\end{equation*}
Combining this with \eqref{X1}-\eqref{X3}, we arrive at
\beg{equation*}\beg{split} \GG(r)&:= \int_0^r \e^{-2\ll t}  \E\big[1_{\{t<\tau_m\}} |X(t)-Y(t)|^2\big]\,\d t\\
&\le \int_0^r \e^{-2\ll t}  \E\Big\{1_{\{t<\tau_m\}}\Big(8|\LL(t)|^2+ \ff 8 7 |\Xi(t)|^2\Big)\Big\}\,\d t\\
&\le \ff 6 7 \GG(r) +\ff 87 C_0\int_0^r\GG(t)\d t + 8 K(T) \int_0^r\GG(t)\d t\\
&\le \ff 6 7 \GG(r) + 8(C_0+K(T)) \int_0^r\GG(t)\d t,\ \ r\in [0,T].\end{split}\end{equation*}
Since by the definitions of $\GG$ and $\tau_m$ we have $\GG(t)<\infty$ for $t\in [0,T]$, it follows from   Gronwall's  inequality that  $\GG(T)=0$. Therefore,   \eqref{3.5} holds and the proof is finished.
\end{proof}

\paragraph{Remark 2.1.} Due to the unbounded term $\ll-A$ in the definition of $\Xi(t)$, even when the time delay is bounded   we are not able to prove Proposition \ref{P3.3} 
with the following weaker  condition in place of $(A4')$:
\beq\label{A4"} |B(t,\xi)-B(t,\eta)|\le C\|\xi-\eta\|_\infty,\ \ t\in [0,T],  \xi,\eta\in C([-r_0,0];\H).\end{equation}
However, when $\H$ is finite-dimesional, $A$ becomes bounded  so that the proof of  Proposition \ref{P3.3} can be modified by using \eqref{A4"} in place of $(A4')$. This will be discussed  in a forthcoming paper. 

\subsection{Proof of Theorem \ref{T2.1}}

  Let $X_0=\xi\in \C_\nu$ be fixed.

  (a) We first assume that \eqref{nu}, $(A1)$ and $(A2')$-$(A4')$ hold.  Consider the following O-U type  SPDE on $\H$:
\beg{equation*}
\d Z^{\xi}(t)= A Z^{\xi}(t)\d t +Q(t,Z^{\xi}(t))\d W(t),\ \ Z^{\xi}(0)=\xi(0).
\end{equation*}
It is classical by \cite{DZ} that   our assumptions imply the existence, uniqueness and non-explosion of the mild solution:
\beg{equation*}
Z^{\xi}(t)= \e^{At} \xi(0) +\int_{0}^{t} \e^{A(t-s)}Q(s,Z^{\xi}(s))\d W(s),\ \ t\ge 0.
\end{equation*}
Letting $Z^\xi_0=\xi$ (i.e. $Z^\xi(\theta)=\xi(\theta)$ for $\theta \le 0$),  and  taking
\beg{equation*}\beg{split}
&W^{\xi}(t)=W(t)-\int_{0}^{t}\psi(s)\d s,\\
&\psi(s)= \big\{Q^{\ast}(QQ^{\ast})^{-1}\big\}(s,Z^{\xi}(s)) \big\{b(s,Z^\xi(s))+ B(s,Z_s^\xi)\big\},\ \ s,t\in[0,T],\end{split}
\end{equation*} we have  \beg{equation*}\begin{split}
&Z^{\xi}(t)= \e^{At} \xi(0)+\int_{0}^{t} \e^{A(t-s)}B(s,Z_{s}^{\xi})\d s\\
&+\int_{0}^{t} \e^{A(t-s)}b(s,Z^{\xi}(s))\d s+\int_{0}^{t} \e^{A(t-s)}Q(s,Z^{\xi}(s))\d W^{\xi}(s), \ \ t\in[0,T].
\end{split}\end{equation*}
By the Girsanov theorem, $W^{\xi}(t)_{t\in[0,T]}$ is a cylindrical Brownian motion on $\mathbb{\bar{H}}$ under probability $\d\mathbb{Q}^{\xi}=R^{\xi}\d \mathbb{P}$, where
\beg{equation*}
R^{\xi}=\exp\bigg[\int_{0}^{T}\big\langle \psi(s),\d W(s)\big\rangle_{\mathbb{\bar{H}}}-\frac{1}{2}\int_{0}^{T}\big|\psi(s)\big|_{\mathbb{\bar{H}}}^{2}\d s\bigg].
\end{equation*}
Then, under the probability $\Q^\xi$,  $(Z^{\xi}(t),W^{\xi}(t))_{t\in[0,T]}$ is a weak mild solution to \eqref{E1}.  On the other hand, by Proposition \ref{P3.3},   the pathwise uniqueness holds for the mild solution to \eqref{E1}. So, by the Yamada-Watanabe principle \cite{TS} (see \cite[Theorem 2]{MO} or \cite{CM} for the result in infinite dimensions), the equation \eqref{E1} has a unique mild solution. Moreover, in this case the solution is non-explosive.

(b) In general, take $\psi\in C_{b}^{\infty}([0,\infty))$ such that $0\leq \psi\leq 1$, $\psi(r)=1$ for $r\in[0,1]$ and $\psi(r)=0$ for $r\in[2,\infty]$. For any $m\geq1, t\geq 0, z\in\mathbb{H}, \xi\in\C_\nu$, let
\beg{equation*}\beg{split}
&b^{[m]}(t,z)=b(t,z)\psi\big(m^{-1}|z|\big),\\
&Q^{[m]}(t,z)= Q\big(t,\psi(m^{-1}|z|)z\big),\\
&B^{[m]}(t,\xi)=B(t,\xi)\psi\big(m^{-1}\|\xi\|_{\C_\nu}\big).\end{split}
\end{equation*}
Then $(A2)$-$(A4)$ imply that   $B^{[m]}, Q^{[m]}, b^{[m]}$ satisfy $(A2')$-$(A4')$. Here, we only verify $(A4')$ for $B^{[m]}$ since the other two conditions are obvious for 
$Q^{[m]}$ and $b^{[m]}$ in place of $Q$ and $b$. For any
$\xi,\eta\in\C_\nu$, let for instance $\|\xi\|_{\C_\nu}\ge \|\eta\|_{\C_\nu}$.  By the choice of $\psi$,
$(A4)$ implies
\beg{equation*}\beg{split} &|B^{[m]}(t, \xi)-B^{[m]}(t,\eta)|\\
&\le \psi\big(m^{-1}\|\xi\|_{\C_\nu}\big)|B(t,\eta)-B(t,\xi)|+|B(t,\eta)|\cdot
 \Big|\psi\big(m^{-1}\|\xi\|_{\C_\nu}\big)- \psi\big(m^{-1}\|\eta\|_{\C_\nu}\big)\Big|\\
 &=  1_{\{\|\xi\|_{\C_\nu}\le 2 m\}}\psi\big(m^{-1}\|\xi\|_{\C_\nu}\big)|B(t,\eta)-B(t,\xi)|\\
 &\quad + 1_{\{\|\eta\|_{\C_\nu}\le 2 m\}}|B(t,\eta)|\cdot
 \Big|\psi\big(m^{-1}\|\xi\|_{\C_\nu}\big)- \psi\big(m^{-1}\|\eta\|_{\C_\nu}\big)\Big|\\
 &\le C(m) \|\xi-\eta\|_{\C_\nu}\end{split}\end{equation*}   for some constant $C(m)>0.$ Thus, $(A4')$ holds for $B^{[m]}$ in place of $B$. 

Now, by (a),  equation \eqref{E1} for $B^{[m]}, Q^{[m]}, b^{[m]}$ in place of $B, Q, b$ has a unique non-explosive mild solution $X^{[m]}(t)$ starting at $X_{0}=\xi.$  Let
\beg{equation*}
\tau_{0}:=0,\ \  \tau_{m}=m\wedge\inf\{t\geq 0:\|X^{[m]}_t\|_{\C_\nu}\geq m\},\ m\geq 1.
\end{equation*}
Since $B^{[m]}(s,\xi)=B(s,\xi)$, $Q^{[m]}(s,\xi(0))=Q(s,\xi(0))$ and $b^{[m]}(s,\xi(0))=b(s,\xi(0))$ hold for $s\leq m$ and $\|\xi\|_{\C_\nu}\leq m$,  Proposition 3.3 implies     $X^{[m]}(t)=X^{[n]}(t)$
for any $n, m\geq1$ and  $t\in[0,\tau_{m}\wedge\tau_{n}]$. In particular, $\tau_{m}$ is increasing in $m$. Let $\zeta=\lim_{m\to\infty}\tau_{m}$ and
\beg{equation*}
X(t)=\sum_{m=1}^{\infty}1_{[\tau_{m-1},\tau_{m})}X^{[m]}(t),\ \ t\in[0,\zeta).
\end{equation*}
It is easy to see that $X(t)_{t\in[0,\zeta)}$ is a mild solution to \eqref{E1} with lifetime $\zeta$, since condition  \eqref{nu} and the definition of $\zeta$ imply $\lim_{t\uparrow\zeta}|X(t)|=\infty$ on $\{\zeta<\infty\}$.  Finally,
by Proposition \ref{P3.3}, the mild solution is unique. Then the proof of  Theorem \ref{T2.1}(1) is finished.

(c) Under the conditions of Theorem \ref{T2.1}(2), for a mild solution   $X(t)_{t\in[0,\zeta)}$   to \eqref{E1} with lifetime $\zeta$,  we intend to prove $\zeta=\infty$ a.s. Obviously,  $$\bar X(t):=\int_{0}^{t} \e^{A(t-s)}Q(s,X(s))\d W(s), t\in[0,\zeta),\ \ t\ge 0$$   is an adapted continuous process on $\mathbb{H}$ up to the lifetime $\zeta$. Let $\bar X(t)=0$ for $ t\in(-\infty,0].$ We see that $Y(t):=X(t)-\bar X(t)$ is a mild solution to the equation
\beg{equation*}
\d Y(t)= (A Y(t)+b(t,Y(t)+\bar X(t))+B(t,Y_{t}+\bar X_{t}))\d t,Y_{0}=X_{0},\ \ t\in[0,\zeta).
\end{equation*}
Due to \eqref{2.1}, the increasing property of $h$ and $\Phi$, and noting that $A\leq 0$, this implies that for any $T>0$,
\beg{equation*}\beg{split}
&\d |Y(t)|^{2}\leq 2\big\langle Y(t), b(t,Y(t)+\bar X(t))+B(t,Y_{t}+\bar X_{t}))\big\rangle\d t\\
&\leq2\big\{\Phi_{\zeta\wedge T}(\|Y_{t}\|_{\C_\nu}^{2})+h_{T}(\|\bar X_{t}\|_{\C_\nu})\big\}\d t,\ \  Y_{0}=X_{0},t\in[0,\zeta\wedge T).
\end{split}\end{equation*}
Then
\beg{equation}\label{XX}\beg{split}
|Y(t)|^{2}\leq& |X(0)|^{2}+2\int_{0}^{t}h_{T}(\|\bar X_{s}\|_{\C_\nu})\d s
+2\int_{0}^{t}\Phi_{T}(\|Y_{s}\|_{\C_\nu}^{2})\d s, \  \ t\in [0, T\land \zeta).
\end{split}\end{equation} Since $Y_0=X_0$,    \eqref{nu} implies
\beg{equation*}\beg{split} \|Y_s\|_{\C_\nu}^2 & =  |Y(s)|^2 + \int_{-\infty}^{-s} |Y(s+\theta)|^2\nu(\d \theta) + \int_{-s}^{0} |Y(s+\theta)|^2\nu(\d \theta)\\
&\le \big\{1+\nu([-s,0))\big\} \sup_{r\in [0,s]} |Y(r)|^2 +\kk(s) \int_{-\infty}^0 |X_0(\theta)|^2 \nu(\d\theta)\\
&\le \kk(T) \|X_0\|_{\C_\nu}^2 + \big\{1+\nu([-T,0)) \big\}\sup_{r\in [0,s]}|Y(r)|^2\\
&=:K_1(T)+K_2(T)\sup_{r\in [0,s]}|Y(r)|^2,\ \ s\in [0,T].\end{split}\end{equation*} 
So, by letting
\beg{equation*}\beg{split} H(t)= \sup_{r\in [0,t]}|Y(r)|^2, \ \ 
\aa(T)= |X(0)|^{2}+2\int_{0}^{T}h_{T}(\|\bar X_{s}\|_{\C_\nu})\d s,\end{split}\end{equation*} we obtain from \eqref{XX} that
$$H(t)\le \aa(T) +2\int_0^t \Phi_T\big(K_1(T)+K_2(T)H(s)\big)\d s,\ \ t\in [0,T\land\zeta).$$
Taking
$$\Psi_T(s)= \int_1^s\ff{\d r}{2\Phi_T(K_1(T)+K_2(T)r)},\ \ s\ge 0,$$ we have $\lim_{s\to\infty} \Psi_T(s)=\infty$ due to the assumption on $\Phi$, so that by Biharis' inequality,
$$  H(t)\le \Psi_T^{-1}(\aa(T)+T)<\infty,\ \ t\in [0,\zeta\land T).$$ Since $\sup_{t\in [0,T\land\zeta)} |\bar X(t)|^2<\infty$ a.s., on the set $\{\zeta\le T\}$ we have
$$\infty = \lim_{t\uparrow\zeta} H(t)\le \Psi_T^{-1}(\aa(T)+T)<\infty,\ \  {\rm a.s.}$$ This means $\P(\zeta\le T)=0$ for all   $T>0$ and hence, $\zeta=\infty$ a.s.

\section{Log-Harnack inequality and gradient estimate}

Throughout of this section, we assume that $\H$ is finite-dimensional and  the length of time delay is finite. Since the log-Harnack inequality implies the strong Feller property (see \cite[Theorem 1.4.1]{WB}), and it is easy to
 see that $P_T$ is strong Feller only if supp$\,\nu\subset [-T,0]$, we see that the restriction on bounded time delay is essential for the study.  On the other hand, although   the restriction on finite-dimensions might be technical rather than necessary, we are not able to drop it in the moment. The reason is that   we adopt the argument of \cite{W11} using coupling by change of measures,
 for which  the Hilbert-Schmidt norm of the diffusion coefficient is used. This reduces the framework to finite-dimensions as the diffusion coefficient in \eqref{6.4} below is merely Lipschitz in the operator norm.  
We remark that for SPDEs with Dini drifts but without delay, the log-Harnack inequality is presented in \cite{W} by  using  finite-dimensional approximation and   It\^o's formula as in \cite{RW, WZ}. However,  in the case with delay   the Markov semigroup is associated to the segment solution, for which the corresponding  It\^o formula is not yet available.

  Let $r_0\in (0,\infty)$ such that supp$\,\nu\subset [-r_0,0]$. In this case $\C_\nu$ is reformulated as
$$\C_\nu=\bigg\{\xi: [-r_0,0]\to \H \ \text{is\ measurable\ with}\ \nu(|\xi|^2):=\int_{-r_0}^0 |\xi(\theta)|^2\nu(\d\theta)<\infty\bigg\}.$$
For $f\in \B_b(\C_\nu)$, the length of the gradient of $f$ at point $\xi\in \C_\nu$ is defined by
$$|\nn f|_{\C_\nu}(\xi)= \limsup_{\eta\to \xi} \ff{|f(\eta)-f(\xi)|}{\|\xi-\eta\|_{\C_\nu}}.$$

\beg{thm}\label{T3.1} Let $\H$ be finite-dimensional and ${\rm supp}\,\nu\subset [-r_0,0]$ for some $r_0\in (0,\infty).$  Assume $\eqref{nu}$, $(A1)$ and $(A2')$-$(A4')$.     Then
 for any $s\ge 0$ there exists a constant $C>0$ such that
 the log-Harnack inequality
\beq\label{3.11}
P_{s,T+s+r_0}\log f(\eta)\leq \log P_{s,T+s+r_0} f(\xi)+\ff C{T\land 1} \|\xi-\eta\|_{\C_\nu}^2, \qquad \xi,\eta\in\C_\nu, T>0
\end{equation}
holds for strictly positive functions $f\in\B_{b}(\C_\nu)$.  Consequently,   
\beq\label{3.2}  |\nn P_{s, s+r_0+T} f|_{\C_\nu}^2 \le \ff{C}{T\land 1}\big\{ P_{s, s+r_0+T} f^2- (P_{s, s+r_0+T} f)^2\big\},\ \ T>0, s\ge 0, f\in \B_b(\C_\nu).\end{equation}
\end{thm}

\beg{proof} According to \cite[Proposition 2.3]{ATW14}, the $L^2$-gradient estimate \eqref{3.2} follows from the  log-Harnack inequality \eqref{3.11}.
Moreover, according to the semigroup property and the Jensen inequality, it suffices to prove the log-Harnack inequality for $T\in (0,1].$ Finally, without loss of generality, we may and do assume that $s=0$.

Now, we use Lemma \ref{L3.1} to transform   \eqref{E1}  into a SDE  with regular coefficients.  Let $\{u(t,\cdot)\}_{t\in [0,T]}$ be in Lemma \ref{L3.1} for fixed $T>0$ and  let $u(\theta,\cdot)=u(0,\cdot)$ for $\theta\in [-r_0,0]$. Define
\beg{equation*}\beg{split} & \Theta(t,x)= x+u(t,x),\  \ t\in [-r_0,T], x\in\H,\\
& (\Theta_t(\xi))(\theta)= \Theta(t+\theta, \xi(\theta)),\ \ \theta\in [-r_0.0], t\in [0,T], \xi\in \C_\nu.\end{split}\end{equation*}
Then  there exists $\ll(T)>0$ such that for any $\ll\ge\ll(T),$ $\Theta(t,\cdot)$ is a diffeomorphism on $\H$ for $t\in (-\infty,T]$ such that
\beq\label{XX1}\|\nn u\|_{T,\infty}\le \ff 1 2,\ \  \|\nn\Theta\|_{T,\infty}+ \|\nn \Theta^{-1}\|_{T,\infty}\le 2, \end{equation} where and in what follows, denote $\Theta^{-1}(t,x)= \{\Theta(t,\cdot)\}^{-1}(x)$.  Obviously, $\Theta_t: \C_\nu\to\C_\nu$ is invertible with
$$\{\Theta_t^{-1}(\xi)\}(\theta)= \Theta^{-1}(t+\theta,\xi(\theta)),\ \ \theta\in [-r_0,0], t\in [0,T].$$
Moreover, for a mild solution   $X^\xi(t)$ to \eqref{E1} with $X_0=\xi$, Lemma \ref{L3.1}(2) implies that $Y^\xi(t):= \Theta(t,X^\xi(t))$ solves the equation
  \beq\label{XX2}\beg{split}
Y^{\xi}(t)=\ & \e^{At} Y^{\xi}(0)+\int_{0}^{t} \e^{A(t-s)}\Big\{\big(\nn\Theta(s,\cdot)\big)Q(s,\cdot)\Big\}\big(\Theta^{-1}(s,Y^\xi(s))\big)  \d W(s)\\
&+\int_{0}^{t}(\lambda-A)\e^{A(t-s)}u \big(s, \Theta^{-1}(s, Y^{\xi}(s))\big)\d s\\
&+\int_{0}^{t}\e^{A(t-s)}\big\{\nn \Theta\big(s,\Theta^{-1}(s,Y^\xi(s))\big)\big\}B(s, \Theta_{s}^{-1}(Y^{\xi}_{s}))\d s,\ \
 t\in[0,T].
\end{split}\end{equation}
So, letting
\beg{equation}\label{XXX}\beg{split}
&\tt{Q}(t,x)=\big\{\big(\nn \Theta(t,\cdot)\big)Q(t,\cdot)\big\}\big(\Theta^{-1}(t, x)\big), \ \   x\in\mathbb{H};\\
&\tt{B}(t,\xi)= A\xi(0) +(\ll-A) u\big(t,\Theta^{-1}(t,\xi(0))\big)\\
&\qquad\qquad + \big\{\nn\Theta\big(t, \Theta^{-1}(t, \xi(0))\big)\big\} B\big(t,\Theta_t^{-1}(\xi)\big),\ \  \xi\in\C_\nu,
\end{split}\end{equation}
we have
$$\d Y^{\xi}(t)= \tt{B}(t,Y^{\xi}_{t})\d t+\tt{Q}(t, Y^{\xi}(t))\d W(t),\ \
 t\in[0,T], Y^{\xi}_{0}=\Theta_0(\xi).$$ Thus, 
  $\tt{X}^{\xi}(t):=Y^{\Theta_{0}^{-1}(\xi)}(t)$ solves the following SDE with delay:
\beq\label{6.4}
\d \tt{X}^{\xi}(t)= \tt{B}(t,\tt{X}^{\xi}_{t})\d t+\tt{Q}(t,\tt{X}^{\xi}(t))\d W(t),\ \ t\in[0,T],\tt{X}^{\xi}_{0}=\xi.
\end{equation} Since   $(A2')$,  $(A4')$ and \eqref{XX1} imply
\beq\label{XXX1}\beg{split} &\|\tt Q\|_{T,\infty} + \|(\tt Q\tt Q^*)^{-1}\|_{T,\infty}\le K,\\
& \|\tt Q(t,x)-\tt Q(t,y)\|\le K|x-y|,\ \ x,y\in \H, t\in [0,T],\\
&|\tt B(t,\xi)-\tt B(t,\eta)|\le K\|\xi-\eta\|_{\C_\nu},\ \ \xi,\eta\in\C_\nu, t\in [0,T]
\end{split}\end{equation}
for some constant $K>0$,  this equation has a unique non-explosive  mild solution for any initial value $\xi\in\C_\nu$.
Let $\tt{P}_{t}f(\xi)=\mathbb{E}f(\tt{X}_{t}^{\xi})$. Since $X^\xi(t)= \Theta^{-1}(t,Y^\xi(t)) =\Theta^{-1}(t,\tt X^{\Theta_0(\xi)}(t)),$ we have
\beg{equation*}\beg{split}
&P_{t}f(\xi):=\mathbb{E}f(X_{t}^{\xi})=\mathbb{E}(f\circ\Theta_{t}^{-1})(Y^{\xi}_{t})=
\mathbb{E}(f\circ\Theta_{t}^{-1})(\tt{X}^{\Theta_{0}(\xi)}_{t})\\
&=\tt{P}_{t}(f\circ\Theta_{t}^{-1})(\Theta_{0}(\xi)), \ \ \xi\in\C_\nu,t\in(0,T],f\in \B_{b}(\C_\nu).
\end{split}\end{equation*}
Therefore, by \eqref{XX1}, the desired log-Harnack inequality for $P_{r_0+T}$ follows from the corresponding inequality for $\tt P_{r_0+T}$, which is ensured by the
following Lemma \ref{LL}. \end{proof}

The following result is parallel to \cite[Theorem 4.3.1]{WB} where the uniform norm on the segment space is used instead of $\|\cdot\|_{\C_\nu}$.

\beg{lem}\label{LL} Let $\tt P_t$ be associated to \eqref{6.4} with coefficients satisfying $\eqref{XXX1}$. Then there exists a constant $C>0$ such that
 the log-Harnack inequality
$$
\tt P_{T+r_0}\log f(\eta)\leq \log \tt P_{T+r_0} f(\xi)+ C\bigg( \ff 1 T |\xi(0)-\eta(0)|^2    + \|\xi -\eta \|^2_{\C_\nu}\bigg)$$ holds for all 
$ \xi,\eta\in\C_\nu, T\in (0,1],$
and  strictly positive functions $f\in\B_{b}(\C_\nu)$.\end{lem}

\beg{proof}  The result can be proved in a similar way as in the proof of  \cite[Theorem 4.3.1]{WB}  using coupling by change of measures, the only difference is to use $\|\cdot\|_{\C_\nu}$ in place of
$\|\cdot\|_\infty$. We remark that the coupling by change of measures are 
  developed in \cite{ATW06} and \cite{W10} to prove  the dimension-free Harnack inequality and   the log-Harnack inequality respectively, see \cite{WB} for more results and discussions.
 For completeness, we include below a brief proof.

(a) For any $\xi,\eta\in \C_\nu$, let $X(t)=\tt X^\xi(t)$ solve \eqref{6.4}, and let $Y(t)$ solve the following SDE with delay:
\beq\label{X0} \beg{split}\d Y(t)= &\tt B(t,X_t)\d t +\tt Q(t,Y(t))\d W(t) \\
&+\ff{1_{[0,T)}(t)}{\gg(t)} \tt Q(t,Y(t))\big\{\tt Q^*(\tt Q\tt Q^*)^{-1}\big\}(t,X(t))\big\{X(t)-Y(t)\big\}\d t,\ \ t\ge 0, Y_0=\eta,\end{split}\end{equation}
where 
$$\gg(t):= \ff 1 {K^2} \big(1- \e^{(t-T)K^2}\big),\ \ t\in [0,T].$$ 
Here, following the line of \cite{ES},  we take the delay term $\tt B(t,X_t)$ instead of $\tt B(t,Y_t)$ such that the SDE for $X(t)-Y(t)$ does not have time delay.    Hence, for any solution $\{Y(t)\}_{t\ge 0}$ to \eqref{X0} with coupling time:
$$\tau:= \inf\{t\in [0,T]:\ X(t)=Y(t)\},\ \ \inf\emptyset:=\infty,$$ the modified process
$$\tt Y(t):= Y(t)1_{\{t<\tau\}} + X(t) 1_{\{t\ge \tau\}}$$ solves \eqref{X0} as well. Using $\tt Y$ to replace $Y$ we may and do assume   $Y(t)=X(t)$ for $t\ge \tau$, so that $X_{T+r_0}=Y_{T+r_0}$  provided $\tau\le T$.  
This is crucial to derive the log-Harnack inequality. Moreover, we take the additional unbounded drift term in \eqref{X0}   to ensure $\tau\le T$, and the idea comes from \cite{W11}.

(b) Since the coefficients in \eqref{X0} are Lipschtiz continuous in the space variable locally uniformly in $t\in [0,T)$, it has a unique solution up to time $T$. To construct a solution for all $t\ge 0$, 
we reformulate the equation as \eqref{6.4} using Girsanov transform. Let
\beq\label{*D00}\beg{split}  \phi(t)=  &\big\{\tt Q^*(\tt Q\tt Q^*)^{-1}\big\}(t,Y(t))\big\{\tt B(t,Y_t)-\tt B(t,X_t)\big\}\\
&-\ff{1}{\gg(t)}\big\{\tt Q^*(\tt Q\tt Q^*)^{-1}\big\}(t,X(t))\big\{X(t)-Y(t)\big\},\ \ t\in [0,T).\end{split}\end{equation}
By \eqref{XXX1} we have 
\beq\label{*D0} |\phi(t)|_{\bar\H}\le K_0\|X_t-Y_t\|_{\C_\nu}+ \ff {K_0|X(t)-Y(t)|}{\gg(t)},\ \ t\in [0,T)\end{equation} for some constant $K_0>0$. Since \eqref{XXX1} implies assumption {\bf (A4.4)} in \cite{WB} for $K_4=K^2$, it follows from \cite[(i) on page 92]{WB} that
\beq\label{*D1} R(t):= \exp\bigg[\int_0^t \<\phi(s), \d W(s)\>_{\bar\H}-\ff 1 2\int_0^t |\phi(s)|^2_{\bar\H} \d s\bigg],\ \ t\in [0,T)\end{equation} is a uniformly integrable martingale such that
\beq\label{*D2} R:=\lim_{t\uparrow T} R(t)\end{equation} exists, and $\d\Q:= R\d\P$ is a probability measure on $\OO$. Moreover, by the Girsanove theorem,
\beq\label{*D0'} \tt W(t):= W(t)-\int_0^{t\land T} \phi(s)\d s,\ \ t\ge 0\end{equation}  is a cylindrical  Brownian motion on $\bar\H$ under probability $\Q$, and according to \cite[(ii) on page 92]{WB} we have $\tau\le T, \Q$-a.s. So, as explained in the end of (a), we have 
\beq\label{X02} X_{T+r_0}=Y_{T+r_0},\ \ \Q\text{-a.s.}\end{equation}
Now, by \eqref{*D00} and \eqref{*D0'} we reformulate \eqref{X0} as
\beq\label{X00} \d Y(t)= \tt B(t, Y_t)\d t +\tt Q(t,Y(t))\d\tt W(t),\ \ Y_0=\eta.\end{equation} By the weak uniqueness of \eqref{6.4}, we have
$\tt P_{T+r_0}(\log f)(\eta)=\E_\Q[\log f(Y_{T+r_0})].$ According to Young inequality (see \cite[Lemma 2.4]{ATW09}), this together with \eqref{X02} implies
\beq\beg{split} \label{X03} \tt P_{T+r_0}(\log f)(\eta)&= \E[R\log f(Y_{T+r_0})]= \E[R\log f(X_{T+r_0})]\\ &\le \log \E f(X_{T+r_0}) + \E[R\log R]
 =\log \tt P_{T+r_0} f(\xi)+ \E_\Q\log R.\end{split}\end{equation}

 (c) To estimate $\E_\Q\log R$,  by \eqref{*D00} and \eqref{*D0'} we reformulate the equation \eqref{6.4} for $X(t)$ as
 \beg{equation*}\beg{split} \d X(t)= & \tt B(t,X_t)\d t +\tt Q(t,X(t))\d\tt W(t)\\
 &- \tt Q(t,X(t))\big\{\tt Q^*(\tt Q\tt Q^*)^{-1}\big\}(t,Y(t))\big\{\tt B(t,X_t)-\tt B(t,Y_t)\big\}\d t\\
 &- \ff {X(t)-Y(t)}{\gg(t)} \d t,\ \ t\in [0,T).\end{split}\end{equation*}
 Combining this with \eqref{X00} we obtain
 \beg{equation*}\beg{split} \d (X(t)-Y(t))= &  \big\{\tt Q(t,X(t))-\tt Q(t, Y(t))\big\}\d\tt W(t)\\
 &+ \Big[I-\tt Q(t,X(t))\big\{\tt Q^*(\tt Q\tt Q^*)^{-1}\big\}(t,Y(t))\Big] \big\{\tt B(t,X_t)-\tt B(t,Y_t)\big\}\d t\\
 &- \ff {X(t)-Y(t)}{\gg(t)} \d t,\ \ t\in [0,T).\end{split}\end{equation*}
 By It\^o's formula and \eqref{XXX1}, there exists a constant $C_0>0$ such that
 \beq\label{M1}\beg{split} \d |X(t)-Y(t)|^2 \le \Big\{& C_0 \|X_t-Y_t\|_{\C_\nu}|X(t)-Y(t)| +K^2|X(t)-Y(t)|^2\\
 &- \ff{2|X(t)-Y(t)|^2}{\gg(t)}\Big\}\d t +\d M(t),\ \ t\in [0,T), \end{split}\end{equation} where
 $\d M(t):= 2 \<(\tt Q(t,X(t))-\tt Q(t, Y(t)))\d \tt W(t), X(t)-Y(t)\>$ is a $\Q$-martingale. Since for some constant $C_1>0$ we have 
\beg{equation*}\beg{split}&C_0 \|X_t-Y_t\|_{\C_\nu}|X(t)-Y(t)|+K^2|X(t)-Y(t)|^2-\ff{2|X(t)-Y(t)|^2}{\gg(t)} \\
& \le C_0 \|X_t-Y_t\|_{\C_\nu}|X(t)-Y(t)|-K^2|X(t)-Y(t)|^2\\
& \le C_1\|X_t-Y_t\|_{\C_\nu}^2,\ \ t\in [0,T),\end{split}\end{equation*}
 \eqref{M1} and  \eqref{nu}  imply 
 \beg{equation*}\beg{split} \E_\Q\|X_t-Y_t\|^2_{\C_\nu} &= \E_\Q|X(t)-Y(t)|^2 +\int_{-r_0}^0 \E_\Q|X(t+\theta)-Y(t+\theta)|^2\nu(\d\theta)\\
 &\le |\xi(0)-\eta(0)|^2 +C_1\int_0^t \E_\Q\|X_s-Y_s\|_{\C_\nu}^2\d s \\
 &\quad+\kk(T) \int_{-r_0}^0 |\xi(\theta)-\eta(\theta)|^2\nu(\d\theta) + \nu([-r_0,0)) \int_0^t \E_\Q|X(s)-Y(s)|^2\d s\\
 &\le C_2 \|\xi-\eta\|_{\C_\nu}^2 +C_2 \int_0^t \E_\Q \|X_s-Y_s\|_{\C_\nu}^2\d s,\ \ t\in [0,T)\end{split}\end{equation*} for some constant $C_2>0$. Therefore, by Gronwall's lemma, 
 \beq\label{M2} \E_\Q\|X_t-Y_t\|_{\C_\nu}^2 \le C_2 \e^{C_2T} \|\xi-\eta\|_{\C_\nu}^2,\ \ t\in [0,T).\end{equation}

 On the other hand, by \eqref{M1} and noting that $2+\gg'-K^2\gg=1$,  we have
 \beg{equation*}\beg{split} &\d \bigg\{\ff{|X(t)-Y(t)|^2}{\gg(t)}\bigg\}\\
 & \le \bigg\{\ff{C_0\|X_t-Y_t\|_{\C_\nu} |X(t)-Y(t)|}{\gg(t)}-\ff{|X(t)-Y(t)|^2}{\gg(t)^2}\big(2+\gg'(t)-K^2\gg(t)\big)\bigg\}\d t+\ff 1{\gg(t)}\d M(t)\\
 &\le  \bigg\{\ff{C_0^2}2 \|X_t-Y_t\|_{\C_\nu}^2 -\ff{|X(t)-Y(t)|^2}{2\gg(t)^2}\bigg\}\d t+\ff 1{\gg(t)}\d M(t),\ \ t\in [0,T).\end{split}\end{equation*}
 So, it follows from \eqref{M2} that
 $$\E_\Q\int_0^T\ff{|X(t)-Y(t)|^2} {\gg(t)^2} \d t\le \ff{2|\xi(0)-\eta(0)|^2}{\gg(0)} + C_0^2 \int_0^T \E_\Q\|X_t-Y_t\|_{\C_\nu}^2\d t.$$ Combining this with \eqref{*D0}, \eqref{*D1}, \eqref{*D2} and \eqref{M2}, we arrive at
\beg{equation*}\beg{split}\E_\Q\log R &= \lim_{t\uparrow T} \E_\Q\bigg\{\int_0^t \<\phi(s),\d\tt W(s)\>_{\bar\H}+\ff 1 2 \int_0^t |\phi(s)|^2_{\bar\H}\d s\bigg\}\\
&=\ff 1 2 \int_0^T \E_\Q|\phi(s)|^2_{\bar\H}\d s \le  C\bigg( \ff 1 T |\xi(0)-\eta(0)|^2   + \ |\xi-\eta|^2_{\C_\nu}  \bigg)\end{split}\end{equation*}
for some constant $C>0.$  Then the proof is finished by \eqref{X03}.

\end{proof}


\beg{thebibliography}{99}

\bibitem{ATW06} 	M.  Arnaudon, A. Thalmaier, F.-Y. Wang,  \emph{Harnack inequality and heat kernel estimate on manifolds with curvature unbounded below,}
 Bull. Sci. Math. 130(2006), 223--233.
 
 \bibitem{ATW09} M. Arnaudon, A. Thalmaier, F.-Y. Wang, \emph{Gradient estimates and Harnack inequalities on non-compact Riemannian manifolds,}
  Stoch. Proc. Appl. 119(2009), 3653--3670.

\bibitem{ATW14} 	M. Arnaudon, A. Thalmaier, F.-Y. Wang, \emph{ Equivalent log-Harnack and gradient for point-wise curvature lower bound,}  Bull. Math. Sci. 138(2014), 643--655.

\bibitem{B} D. Bakry, I. Gentil, M. Ledoux, \emph{Analysis and Geometry of Markov Diffusion Operators,}   Springer,  2014.

\bibitem{DZ} G. Da Prato, J. Zabczyk,  \emph{Stochastic Equations in Infinite Dimensions,} Cambridge University
Press,   1992.

\bibitem{ES} A. Es-Sarhir, M.-K. v. Renesse, M. Scheutzow, \emph{Harnack inequality for functional SDEs with
bounded memory,}  Electron. Commun. Probab. 14(2009), 560--565. 

\bibitem{MO} M. Ondrejet,  \emph{Uniqueness for stochastic evolution equations in Banach spaces,} Dissertationes
Math. (Rozprawy Mat.) 426(2004).
\bibitem{CM} C. Prevot, M. R\"ockner,  \emph{A Concise Course on Stochastic Partial Differential Equations,} Lecture Notes in Math. Vol. 1905, Springer, Berlin, 2007.

\bibitem{PW} E. Priola, F.-Y. Wang, \emph{Gradient estimates for diffusion semigroups with singular coefficients,}  J. Funct. Anal. 236(2006), 244--264.

\bibitem{RW} M. R\"ockner, F.-Y. Wang, \emph{Log-Harnack  Inequality for Stochastic differential equations in Hilbert spaces and its consequences, } Infinite Dimensional Analysis, Quantum Probability and Related Topics 13(2010), 27--37.

\bibitem{SWY} J. Shao, F.-Y. Wang, C. Yuan, \emph{Harnack  inequalities for stochastic (functional) differential equations with non-Lipschitzian coefficients, } Electron. J. Probab.  17(2012), 1--18.

\bibitem{W10} 	F.-Y. Wang, \emph{Harnack inequalities on manifolds with boundary and applications,} J. Math. Pures Appl. 94(2010), 304--321.

\bibitem{W11} 	F.-Y. Wang, \emph{Harnack inequality for SDE with multiplicative noise and extension to Neumann semigroup on non-convex manifolds,} Ann.  Probab. 39(2011), 1449--1467.

\bibitem{WB} F.-Y. Wang, \emph{Harnack Inequality and Applications for Stochastic Partial Differential Equations,} Springer,   2013.

\bibitem{Wbook} F.-Y.  Wang, \emph{Analysis for Diffusion Processes on Riemannian Manifolds, } World Scientific, 2014. 

\bibitem{W} F.-Y. Wang,  \emph{Gradient estimate and applications for SDEs in Hilbert space with multiplicative noise and Dini continuous drift,}  arXiv:1404.2990.

\bibitem{WZ} F.-Y. Wang, T. Zhang,  \emph{  Gradient estimates for stochastic evolution equations with non-Lipschitz coefficients,}  J. Math. Anal. Appl. 365(2010), 1--11.
\bibitem{TS} T.Yamada, S.Watanabe,  \emph{On the uniqueness of solutions of stochastic differential equations,} J. Math. Kyoto Univ. 11(1971), 155--167.

\end{thebibliography}

\end{document}